\documentclass[12pt]{article}
\font\sdopp=msbm10
\def\CI {\sdopp {\hbox{C}}}
\title{
Clifton-Pohl torus and 
geodesic completeness by  
a 'complex' point of view
}
\author{
Claudio Meneghini 
\footnote{\scriptsize
Math.Dept.University of Parma,
Str.M.D'Azeglio 85, 43100 Parma (Italy)
}
}
\begin{document}
\maketitle
\bibliographystyle{plain} 
\parindent=8pt
\font\cir=wncyb10
\def\Iu{\cir\hbox{YU}}
\def\Ze{\cir\hbox{Z}}
\def\pe{\cir\hbox{P}}
\def\Ef{\cir\hbox{F}}
\def\ced{{\hbox{\tt\large\c c}}}
\font\sdopp=msbm10
\def\DI{\sdopp{\hbox{D}}}
\def\ESSE {\sdopp {\hbox{S}}}
\def\ERRE {\sdopp {\hbox{R}}}
\def\CI {\sdopp {\hbox{C}}}
\def\ENNE{\sdopp {\hbox{N}}}
\def\ZETA{\sdopp {\hbox{Z}}}
\def\PI {\sdopp {\hbox{P}}}
\def\M{\hbox{\boldmath{}$M$\unboldmath}} 
\def\N{\hbox{\boldmath{}$N$\unboldmath}} 
\def\P{\hbox{\boldmath{}$P$\unboldmath}} 
\def\Y{\hbox{\boldmath{}$Y$\unboldmath}} 
\def\tr{\hbox{\boldmath{}$tr$\unboldmath}} 
\def\f{\hbox{\boldmath{}$f$\unboldmath}} 
\def\u{\hbox{\boldmath{}$u$\unboldmath}} 
\def\v{\hbox{\boldmath{}$v$\unboldmath}} 
\def\U{\hbox{\boldmath{}$U$\unboldmath}} 
\def\V{\hbox{\boldmath{}$V$\unboldmath}} 
\def\W{\hbox{\boldmath{}$W$\unboldmath}} 
\def\id{\hbox{\boldmath{}$id$\unboldmath}} 
\def\alph{\hbox{\boldmath{}$\aleph$\unboldmath}} 
\def\bet{\hbox{\boldmath{}$\beta$\unboldmath}} 
\def\gam{\hbox{\boldmath{}$\gamma$\unboldmath}} 
\def\U{\mathop{u}\limits}
\def\f{\hbox{\boldmath{}$f$\unboldmath}} 
\def\g{\hbox{\boldmath{}$g$\unboldmath}} 
\def\h{\hbox{\boldmath{}$h$\unboldmath}} 
\def\IM{\hbox{\boldmath{}$i$\unboldmath}} 

\def\Ch{\hbox{\rm Ch}}
\def\Th{\hbox{\rm Th}}
\def\Sh{\hbox{\rm Sh}}
\def\SettTh{\hbox{\rm SettTh}}
\def\CIRC{\mathop{\tt o}\limits}
\def\quan{\vrule height6pt width6pt depth0pt}
\def\QUAN{\ \quan}
\def\BETA{\mathop{\beta}\limits}
\def\GAMMA{\mathop{\gamma}\limits}
\def\VI{\mathop{v}\limits}
\def\UI{\mathop{u}\limits}

\def\oi{\mathop{\omega}\limits}
\def\ei{\mathop{\eta}\limits}
\def\FI{\mathop{\varphi}\limits}
\def\poi{\mathop{\psi}\limits}

\def\VII{\mathop{V}\limits}
\def\WI{\mathop{w}\limits}
\def\ZETA{\mathop{Z}\limits}
\def\BETA{\mathop{\beta}\limits}
\def\GAMMA{\mathop{\gamma}\limits}
\def\VI{\mathop{v}\limits}
\def\UI{\mathop{u}\limits}
\def\VII{\mathop{V}\limits}
\def\WI{\mathop{w}\limits}
\def\ZETA{\mathop{Z}\limits}
\def\ssqrt#1{\left(#1\right)^{1/2}}
\def\sssqrt#1{\left(#1\right)^{-1/2}}
\def\TTT{\sl}
\def\BBB{\sl}

\newtheorem{definition}{Definition}
\newtheorem{defi}[definition]{D\'efinition}
\newtheorem{lemma}[definition]{Lemma}
\newtheorem{lemme}[definition]{Lemme}
\newtheorem{proposition}[definition]{Proposition}
\newtheorem{theorem}[definition]{Theorem}        
\newtheorem{theoreme}[definition]{Th\'eor\`eme}        
\newtheorem{corollary}[definition]{Corollary}  
\newtheorem{corollaire}[definition]{Corollaire}  
\newtheorem{remark}[definition]{Remark}  
\newtheorem{remarque}[definition]{Remarque}
  
\font\sdopp=msbm10
\def\ERRE {\sdopp {\hbox{R}}}
\def\QU {\sdopp {\hbox{Q}}}
\def\CI {\sdopp {\hbox{C}}}
\def\DI {\sdopp {\hbox{D}}}
\def\ENNE{\sdopp {\hbox{N}}}
\def\ZETA{\sdopp {\hbox{Z}}}
\def\PI {\sdopp {\hbox{P}}}

\def\M{\hbox{\tt\large M}}
\def\N{\hbox{\tt\large N}}
\def\T{\hbox{\tt\large T}}

\def\P{\hbox{\boldmath{}$P$\unboldmath}} 
\def\tr{\hbox{\boldmath{}$tr$\unboldmath}} 
\def\f{\hbox{\large\tt f}} 
\def\g{\hbox{\tt\large g}}

\def\F{\hbox{\boldmath{}$F$\unboldmath}} 
\def\G{\hbox{\boldmath{}$G$\unboldmath}} 
\def\L{\hbox{\boldmath{}$L$\unboldmath}} 
\def\h{\hbox{\boldmath{}$h$\unboldmath}} 
\def\e{\hbox{\boldmath{}$e$\unboldmath}} 

\def\u{\hbox{\boldmath{}$u$\unboldmath}} 
\def\v{\hbox{\boldmath{}$v$\unboldmath}} 
\def\U{\hbox{\boldmath{}$U$\unboldmath}} 
\def\V{\hbox{\boldmath{}$V$\unboldmath}} 
\def\W{\hbox{\boldmath{}$W$\unboldmath}} 
\def\id{\hbox{\boldmath{}$id$\unboldmath}} 
\def\alph{\hbox{\boldmath{}$\alpha$\unboldmath}} 
\def\bet{\hbox{\boldmath{}$\beta$\unboldmath}} 
\def\gam{\hbox{\boldmath{}$\gamma$\unboldmath}} 
\def\pphi{\hbox{\boldmath{}$\varphi$\unboldmath}} 
\def\ppsi{\hbox{\boldmath{}$\psi$\unboldmath}} 
\def\Ppsi{\hbox{\boldmath{}$\Psi$\unboldmath}} 
\def\labelle #1{\label{#1}}
\begin{abstract}
We 
show that a natural 
complexification and a mild generalization of
the idea of completeness
guarantee
geodesic completeness of
Clifton-Pohl torus; we explicitely 
compute all of its geodesics.
\end{abstract}

{{\tt Keywords:} Clifton-Pohl torus, geodesic completeness,
holomorphic metric, analytical continuation,
elliptic functions.}\\

\section{Foreword}

The geodesic equations of 
of the so-called
'Clifton-Pohl
torus' $\T$ (a compact, geodesically incomplete, Lorentz manifold, see 
\cite{oneill}, 7.16), can be naturally thought of
as a system of ordinary differential equations
in the complex domain.
Thus they yield in fact holomorphic germs of
solutions: it is therefore a natural idea to
'complexify' the environment to make all
this sound.

Therefore, we propose a natural, 
 and mild, generalization 
of the idea of
geodesic completeness,
amounting to 
conjugating analytical continuation
with complexification of the environment.
We shall 
show that all geodesics of 
$\T$ 
are 'complete' by this point of view: more 
precisely, their germs of solution
will be shown to admit 'endless' 
analytical continuability on $\CI$.

We shall use the concept
of a {\it holomorphic metric}
 on a complex manifold $\M$ 
(see e.g. \cite{lebrun}
): it amounts to a nondegenerating symmetric section 
of
the twice covariant holomorphic tensor bundle 
 $\hbox{\large\sf T}_{0}^{2}\M$.
Of course, it carries no 'signature';
however, by simmetry, it induces
a canonical Levi-Civita's connexion on $\M$, allowing
geodesics to be defined as 
auto-parallel paths; in this paper
we shall not fathom
these issues beyond: the reader 
is referred also to \cite{manin}.
Finally, if $\M$ arises as a 'complexification' of a semi-Riemannian manifold $\N$,
it is easily seen that the real geodesics  of $\N$ are restrictions to the real axis of the complex ones of $\M$ and vice versa (see \cite{lebrun}).
This fact sometimes allows us to 'flank' isolated singularities on the real line by running along complex trips, i.e. to 'connect'
geodesics which, in the usual sense,
are completely unrelated.
In other cases, geodesics which end at finite time by the usual point of view, could be continued only by admitting 'complexified' values.

We suggest an idea
of 
our notion of {\sl completeness}
(see also definition 
\ref{completessa})~: a compex curve $\gamma$ 
from a Riemann
surface $(S,\pi)$ over $\CI$ into
a complex manifold $\M$ will be told to be 
{\bf complex-complete} if
$\CI\setminus\pi(S)$ is a discrete set.
In other words, up to 
changing branch, $\gamma$ is 
analytically continuable everywhere,
except at most at a discrete set in the
complex plane.

{\it En passant}, our procedure will yield
explicit solutions for all geodesics of $ \T  $,
mainly by means of elliptic functions.

\section{Basic definitions and lemmata}
In the following,  ${\cal U}$ 
will be a region in the 
complex plane and $\M$
a complex manifold.
\begin{definition}
{\TTT A 
holomorphic 
metric on} $\M$ is an everywhere maximum-rank 
symmetric section 
of
the twice covariant holomorphic tensor bundle 
 $\hbox{\large\sf T}_{0}^{2}\M$.
A {\BBB holomorphic 
Riemannian manifold} is a complex manifold endowed with a 
holomorphic 
metric.
\label{riemann}
\end{definition}

The idea of the analytical continuation of a holomorphic mapping element
(or of a germ) 
$f:{\cal U}\rightarrow\M$  is well known and amounts to a quadruple $Q_{\M}=(S,\pi,j,F)$, where
$S$ is a connected Riemann surface over a region of $\CI$,
$\pi\,\colon\, S\rightarrow \CI$ is a nonconstant holomorphic mapping 
such that $U\subset \pi(S)$,
$j\,\colon\, U\rightarrow S$ is a holomorphic  immersion such that $\pi\circ j=id\vert_{U}$
and 
$F\,\colon\, S\rightarrow \M$ is a holomorphic mapping such that $F\circ j=f$.
We do not consider 
branch points here;
it is a well known
result that there exists a unique maximal analytical continuation,
called the (regular)
{\BBB Riemann surface}, 
of $\left({\cal U},f   \right)$, which is made up
by all elements which are analytical 
continuations
of $f   $.
As general references, see e.g. \cite{ahlfors}, \cite{cassa},\cite{malgrange},
\cite{narasimhan},
\cite{palka}.
 
\begin{definition}
\label{completezza}
A complex curve $F:S\to\M   $ defined on a 
Riemann surface $\displaystyle
\left(S,\pi,j,F
\right)$ over $\CI$
is {\bf 
complex-complete}
 provided that 
$\CI\setminus
\displaystyle\pi\left(S
    \right)$ is a 
discrete set in the 
complex plane;
a real-analytic Lorentz manifold $ \N  $
is (weakly) complex-complete provided that,
for every geodesic $ \gamma  $ in $ \N  $,
there exists a complexification $ \M_{\gamma}$
of $ \N  $, with embedding mapping 
$\ced_{\gamma}: \N\to\M_{\gamma}$,
such that the Riemann surface of 
$ \ced\circ\gamma  $ is complex-complete.
\labelle{completessa}
\end{definition}

We conclude this section by
remarking that,
as a consequence of the 
existence-and-uniqueness 
theorem of o.d.e.'s theory 
in the complex domain (see e.g. \cite{hille},
 th 2.2.2, \cite{ince} 
p.281-284),
 for each point $p$ in
a holomorphic Riemannian 
manifold and each holomorphic 
tangent vector $X$ at $p$, 
there exists a unique 
holomorphic geodesic element 
starting at $p$ with velocity 
$X$.
 
\section{Clifton-Pohl torus}
Let now 
$\N:=\ERRE^2
\setminus \{0\}$ be endowed
 with the Lorentz metric
${du\odot dv}/({u^2+v^2})
$;
the group $D$ generated by scalar
multiplication by $2$ is a group of isometries of 
$\N$; 
its action is properly dicontinuous, hence
$\T=\N/D$ is a Lorentz surface.
Now, $\T$ is topologically
equivalent to 
the closed annulus 
$1\leq\varrho\leq 2$, 
with boundaries identified 
by the action of $D$, 
i.e. a torus; notwithstanding, $\T$ is geodesically incomplete, since
$t\mapsto\left(1/(1-t),0\right)$
 is a geodesic of 
$\hbox{\tt N}$ (see
\cite{oneill}).
This example suggest us
(in fact
rather weakly)
to extend the domain of the definition
of the natural parameter $t$ to
the complex plane: of course, this in turn requires
a complexification of $\T   $ or, equivalently,
 $\N   $.
Doing this would render 
the above curve complex-complete,
 according to definition 
\ref{completessa}.
In the following, we shall indeed
study the
 holomorphic Riemannian 
manifold
$
\M=\left[\CI^2
\setminus ((1,i)\CI \cup (1,-i)\CI),
{du\odot dv}/({u^2+v^2})
   \right]
$, which is a natural complexification
of $\N$: by methonymy, we shall 
use the name 'Clifton-Pohl
torus'
for
$\N$ and $\M$ too, 
since our completeness 
theorems
can be 
easily be
'restricted' to the real slice
 $\N$ and then pushed down
from $\N$ to $\T$
 with respect to the action of
$D$.
Here is our main result.
\begin{theorem}
Each geodesic 
$\gamma$ of $ \M  $ 
is complex-complete.
\labelle{principal}
\end{theorem}
{\bf Proof.}
The geodesic equations of 
$\M$ are:
\begin{equation}
\UI^{\bullet\bullet}=
{2u}/({u^2+v^2})
\UI^{\bullet}{}^2,\quad 
\VI^{\bullet\bullet}=
{2v}/({u^2+v^2})
\VI^{\bullet}{}^2.
\labelle{equazions}
\end{equation}

Let us first study 
null geodesics of $\M$: 
without loss 
of generality,
it is enough to deal 
with the case
$v\equiv A$.
Equations (\ref{equazions})
imply
$\UI^{\bullet\bullet}=
{2u}/({u^2+A^2})
\UI^{\bullet}{}^2$, which is 
solved by 
$t\mapsto (C-Bt)^{-1}$ if $A=0$
 and by 
$t\mapsto\tan(At+B)$ 
if $A\not=0$, for suitable 
complex constants $B$ and $C$.
The above functions are 
meromorphic, hence complex-complete
by definition
 \ref{completessa}.
We turn to nonnull geodesics:
%
%
let $\gamma$ start
at $(\alpha,\beta)$, 
with velocity
$(x,y)$; 
we may suppose, without loss of generality,
$\alpha\not=0$ and
$\beta\not=0$.
Moreover, we have $x\not=0$ and
$y\not=0$ otherwise $\gamma$ 
would be null.
The equations (\ref{equazions})
can be integrated once to yield:
\begin{equation}
\UI^{\bullet}\VI^{\bullet}=
A(u^2+v^2),\quad u/\UI^{\bullet}+v/\VI^{\bullet}=B
,
\labelle{intprimm}
\end{equation}
where $A=xy/(\alpha^2+\beta^2)$
and $B=\alpha/x+\beta/y$.\\

Now we single out all geodesics
such that $u\VI^{\bullet}\equiv v\UI^{\bullet}$:
by also keeping into account
the equations 
(\ref{intprimm}), they are
easily seen to be of the form
$t\mapsto(a_1\exp(bt),a_2
\exp(bt)) $ for suitable 
complex constants $a_1, a_2$ and
$b$: these curves are clearly 
complex-complete.

Therefore, from now on, 
we may suppose $\alpha y\not=
\beta x$, without 
loss of generality: this implies, by easy 
calculations,
\begin{equation}
AB^2\Ch\log(\alpha/\beta)\not=2,
\labelle{easy}
\end{equation}
for any branch of the logarithm.

Now let us choose a branch \hbox{\tt log}
of the logarithm, defined on a connected neighbourhooud of both $\alpha$ and $\beta$
and
perform
the 
local change of coordinates 
${\omega}=\hbox{\tt log}\, u$,
${\eta}=\hbox{\tt log}\, v$:
the equations
(\ref{intprimm})
are turned into
\begin{equation}
\oi^{\bullet}\ei^{\bullet}=
2A\, \Ch(\omega-\eta),\quad
1/\oi^{\bullet}+1/\ei^{\bullet}=B
.
\labelle{intnoeuv}
\end{equation}
We can solve (\ref{intnoeuv}) with respect to 
$\oi^{\bullet}$ and $\ei^{\bullet}$, getting
\begin{equation}
\cases
{
\oi^{\bullet}=
2\left(B-\sqrt{B^2-2/[A\, \Ch(\oi-\ei)]}
   \right)^{-1}\cr
\ei^{\bullet}=
2\left(B+\sqrt{B^2-2/[A\, \Ch(\oi-\ei)]}
   \right)^{-1}
}.
\labelle{intnoeuv2}
\end{equation}
Subtract and set $\varphi
:=\oi-\ei$; this yields 
\begin{equation}
\FI^{\bullet}=2\,
\sqrt{(AB\Ch\varphi)^2-2A\Ch\varphi}
\labelle{equazion}
\end{equation} 
with the initial 
value $\varphi(0) =
\log(\alpha /\beta )$:
by (\ref{easy}), the differential
equation (\ref{equazion}) is regular at
$\varphi(0)$, and we can choose a branch of
$\SettTh $ at $ \phi_0  $ and
make the further
change of variable $\varphi=2\,\SettTh\,\psi$,
getting, by easy calculations,
\begin{equation}
\displaystyle
\poi^{\bullet}
=
\sqrt{A^2B^2-2A}
\sqrt{(1+\psi^2)\left(1-\frac{A^2B^2+2A}
{2A-A^2B^2}\psi^2\right)}
\labelle{darrera}
\end{equation}
with the (nonsingular) 
initial condition $\psi(0)= \psi_0:=
2\,\SettTh\,\varphi(0)  $.

By separation of variables we can bring back
this equation to the calculation of an elliptic integral, depending on the parameters $A$ and
$B$: as it is shown e.g. in \cite{ahlfors}
p.240, this implies that $ \psi  $ is 
an 
elliptic function $ \hbox{\tt\large e}_{A,B} $,
so that $ \varphi  $ is analytically
continuable everywhere, except at most at the
discrete set $  \hbox{\tt\large e}_{A,B} ^{-1}
(\{ {\pm 1}   \} )$.
Therefore, $ \varphi  $ is complex-complete.

From (\ref{intnoeuv2}),
we get now
that both
$\ei^{\bullet}$ and $\oi^{\bullet}$ is complex-complete; since
integration
clearly preserves completeness,
also $\ei$ and $\oi$ are complex-complete: 
but $u=\exp({\oi})$ and
$v=\exp({\ei})$: this eventually 
implies, by analytical
continuation, that
$\gamma$ is complex-complete.
\QUAN
\noindent
\subparagraph{\bf Final Remarks:}
the concern naturally arises whether there 
exists a minimal complexification of $ \T  $, or 
rather of $ \N  $, which guarantees geodesic completeness: this seems to be a harder
question.

{\tt The author's e-mail 
address: clamengh@
bluemail.ch
}


\begin{thebibliography}{7}
\bibitem[AHL]{ahlfors}
Lars V.Ahlfors,
{\it 'Complex Analysis'}
{McGraw-Hill, 1953}
%

\bibitem[CAS]{cassa}
%
Antonio Cassa,
{\it 'Teoria delle curve algebriche piane e delle superfici di Riemann
compatte'}
{Pitagora,1983}
%
\bibitem[HIL]{hille}
Einar Hille,
{\it 'Ordinary differential equations in the complex domain',}
{John Wiley \& sons, 1976}
%
\bibitem[INC]{ince}
E.L.Ince,
{'Ordinary differential equations'}
{Dover,1956 (originally published in 1926)}
%
\bibitem[LEB]{lebrun}
Claude Lebrun,
{\it 'Spaces of complex null geodesics in complex-Riemannian geometry',}
{Trans. of the AMS, vol 278 n.1, July 1983}
%
\bibitem[MAL]{malgrange}
{B.Malgrange,}
{\it 'Lectures on the theory of functions of several complex variables'}
{Tata institute of fundamental research, Bombay 1958}
%
\bibitem[MAN]{manin}
Yuri Manin,
{\it 'Gauge fields theory and complex geometry'}
{Springer Verlag, 1984}
%
\bibitem[NAR]{narasimhan}
{R. Narasimhan,}
{\it 'Several complex variables'}
{The university of Chicago Press, Chigago and London, 1971}
%
\bibitem[ONE]{oneill}
Barret O'Neill,
{\it 'Semi-Riemannian geometry',}
{Academic Press, 1983}
%
\bibitem[PAL]{palka}
Bruce P.Palka,
{\it 'An Introduction to Complex Function 
Theory'}
{Springer,1991}

\end{thebibliography}
\end{document}